\documentclass[11pt]{amsart}
\usepackage[utf8]{inputenc}
\usepackage[T1]{fontenc}
\usepackage{lmodern}
\usepackage{amsmath,amssymb,amsthm,mathtools}
\usepackage{booktabs}
\usepackage{graphicx}
\usepackage{enumitem}
\usepackage[margin=2.7cm]{geometry}
\usepackage[colorlinks=true,linkcolor=blue,citecolor=blue,urlcolor=blue]{hyperref}

\newtheorem{theorem}{Theorem}[section]
\newtheorem{proposition}[theorem]{Proposition}
\newtheorem{lemma}[theorem]{Lemma}
\newtheorem{corollary}[theorem]{Corollary}
\theoremstyle{definition}
\newtheorem{problem}[theorem]{Problem}
\theoremstyle{remark}
\newtheorem{remark}[theorem]{Remark}

\numberwithin{equation}{section}

\newcommand{\R}{\mathbb{R}}
\newcommand{\dd}{\,\mathrm{d}}
\newcommand{\Hm}{\mathcal{H}^{N-1}}
\DeclareMathOperator{\dv}{div}

\title[Talenti-type rigidity for the Schiffer problem]{A Talenti-type rigidity theorem and P-function constraints for the Pompeiu--Schiffer problem in every dimension}

\author{Diaraf Seck}
\address{Laboratoire de Math\'ematiques de la D\'ecision et d'Analyse Num\'erique (LMDAN), \'Ecole Doctorale de Math\'ematiques et Informatique, Universit\'e Cheikh Anta Diop, BP 16889, Dakar Fann, S\'en\'egal}
\email{diaraf.seck@ucad.edu.sn}

\date{\today}

\begin{document}


\begin{abstract}
Schiffer's conjecture asserts that the ball is the only bounded domain, with boundary homeomorphic to a sphere, carrying a nonconstant Neumann eigenfunction of the Laplacian that is constant on the boundary; for such domains it is equivalent to the Pompeiu problem. In 2026 it was disproved in the plane and in many higher dimensions by computer-assisted and analytic constructions. We study what survives. Our main result is a Talenti-type theorem valid in every dimension $N\ge2$: if the eigenfunction $u$ stays strictly below its boundary value $c>0$, then $kR^*\ge j_{N/2,1}$, where $k^2$ is the eigenvalue, $R^*$ the radius of the ball of equal volume and $j_{N/2,1}$ the first positive zero of $J_{N/2}$; equality holds only for the ball. We compare this threshold with the unconditional bound $kR^*\ge 2^{1/N}j_{N/2-1,1}$ given by the Krahn--Szeg\H{o} inequality. Using Payne's P-function $P=|\nabla u|^2+\lambda u^2$ we prove that every Schiffer domain satisfies $\max|u|>\sqrt N\,|c|$, together with pointwise and distance bounds. We correct the shape-optimization approach of an earlier work by computing the optimality condition of the associated boundary functional with an adjoint state. Using the transcendence of the zeros of Bessel functions of integer and half-integer order, we prove the infinitesimal rigidity of the ball, in every dimension, both for the Schiffer problem and for the Berenstein problem. Finally, we reconstruct numerically the convex four-dimensional counterexample of Guo and check the above statements on it.
\end{abstract}

\maketitle

\section{Introduction}

Let $\Omega\subset\R^N$, $N\ge2$, be a bounded domain whose boundary is smooth and homeomorphic to $S^{N-1}$. Consider the overdetermined problem
\begin{equation}\label{eq:schiffer}
\Delta u+\lambda u=0 \ \text{in } \Omega,\qquad \partial_\nu u=0 \ \text{and}\ u=c \ \text{on } \partial\Omega,\qquad \lambda=k^2>0,\quad u\not\equiv\text{const}.
\end{equation}
We call $\Omega$ a \emph{Schiffer domain} if \eqref{eq:schiffer} has a solution. Schiffer's conjecture asserted that Schiffer domains are balls. The Pompeiu problem \cite{Pompeiu} asks whether a continuous function is determined by its integrals over all rigid motions of $\Omega$; the ball fails this property. Brown, Schreiber and Taylor \cite{BST} and Williams \cite{Williams76} showed that, for such domains, failure of the Pompeiu property is equivalent to the existence of a solution of \eqref{eq:schiffer}: Green's formula then shows that the Fourier transform of $\mathbf 1_\Omega$ vanishes on the sphere $|\xi|=k$. This is Problem~80 in Yau's list \cite{Yau}.

\subsection*{Known results} Berenstein and Yang \cite{BerensteinYang87} proved that infinitely many admissible eigenvalues force the ball, and the boundary of a Schiffer domain is analytic \cite{KinderlehrerNirenberg,Williams81}. In the plane the conjecture is known when $\lambda$ is among the seven lowest Neumann eigenvalues \cite{Aviles,Deng}, when $u$ has no saddle point \cite{WillmsGladwell}, for quadrature domains \cite{Ebenfelt}, for some centrally symmetric convex domains \cite{Mondal,DaiSunWeiZhang}, and the regular pentagon is not a Schiffer domain \cite{NigamSiudejaYoung}. In every dimension, the conjecture holds when the third normal derivative of $u$ is constant on the boundary \cite{Liu}, and for small perturbations of balls \cite{Agranovsky,Kobayashi,Canuto}. Counterexamples were known in other geometries \cite{FallMinlendWeth,CaoLaboraFernandez}.

\subsection*{The 2026 counterexamples} Colbrook and Stepaniants \cite{ColbrookStepaniants} constructed, by a computer-assisted proof, a bounded simply connected planar domain with analytic boundary and $D_{10}$ symmetry which is not a disc and carries a solution of \eqref{eq:schiffer}; an independent family of counterexamples is due to Cao-Labora and de Dios Pont \cite{CaoLaboraDeDiosPont}. Guo \cite{Guo} constructed convex counterexamples in higher dimensions. In the same spirit, Colbrook, Sadeghi and Stepaniants \cite{ColbrookSadeghiStepaniants} disproved the planar Berenstein conjecture: there is a non-circular planar domain carrying a solution of
\begin{equation}\label{eq:berenstein}
\Delta u+k^2u=0 \ \text{in } \Omega,\qquad u=0 \ \text{and}\ \partial_\nu u=1 \ \text{on } \partial\Omega.
\end{equation}

\subsection*{Contributions} The link between these conjectures and shape optimization goes back in particular to \cite{Grigelionis,ChatelainHenrot}. In \cite{Seck24} the author studied them through shape functionals; the present paper re-examines that approach after the 2026 counterexamples and isolates the rigidity statements that survive. Throughout, $R^*=(|\Omega|/\omega_N)^{1/N}$ is the radius of the ball of equal volume, $\omega_N=|B_1|$, and $j_{\mu,n}$ is the $n$-th positive zero of $J_\mu$.

\begin{enumerate}[label=(\roman*)]
\item \emph{A Talenti-type theorem in every dimension} (Theorem~\ref{thm:talenti}). If $c>0$ and $u<c$ in $\Omega$, then $kR^*\ge j_{N/2,1}$, with equality if and only if $\Omega$ is a ball. Without the sign condition, the Krahn--Szeg\H{o} inequality only gives $kR^*\ge 2^{1/N}j_{N/2-1,1}$ (Proposition~\ref{prop:krahn}); see Table~\ref{tab:thresholds}.
\item \emph{Constraints from the P-function} (Section~\ref{sec:P}). Every Schiffer domain satisfies $\max_{\overline\Omega}|u|>\sqrt N\,|c|$ and $|\nabla u|^2+\lambda u^2\le\lambda\max u^2$; the eigenfunction has a critical point where $|u|>|c|$; in the plane, a non-convex Schiffer domain forces a critical point where $|u|<|c|$.
\item \emph{Shape derivatives} (Section~\ref{sec:shape}). The vanishing of the material derivative on the boundary, assumed in \cite{Seck24}, never holds; we compute the correct optimality condition, with an adjoint state, for the boundary functional of \cite{Seck24}.
\item \emph{Infinitesimal rigidity of the ball} (Section~\ref{sec:infinitesimal}), in every dimension, for problem \eqref{eq:schiffer} and for problem \eqref{eq:berenstein}, through the transcendence of the zeros of Bessel functions of integer and half-integer order.
\item \emph{A numerical illustration} (Section~\ref{sec:numerics}): an independent reconstruction of Guo's convex four-dimensional domain, on which the statements above are checked.
\end{enumerate}

The common thread is the sign condition $u<c$: rigidity can be proved as long as a radial test function keeps a constant sign, whereas the known counterexamples violate this condition and live at high frequency.

\section{Preliminaries}\label{sec:prelim}

Let $u$ be a solution of \eqref{eq:schiffer}. Then $c\ne0$, by uniqueness for the Cauchy problem, and
\begin{equation}\label{eq:gradzero}
\nabla u=0\quad\text{on }\partial\Omega .
\end{equation}
In particular $\partial\Omega$ is contained in the critical set of $u$. We denote by $s$ the signed distance to $\partial\Omega$ (negative in $\Omega$, positive outside), by $\nu$ the outward normal and by $H=\dv\nu=\Delta s$ the mean curvature, i.e.\ the sum of the principal curvatures $\kappa_1,\dots,\kappa_{N-1}$ \cite{DelfourZolesio}; we set $|I\!I|^2=\sum\kappa_i^2$, and in dimension $2$ we write $\kappa$ for the curvature. In a neighbourhood of $\partial\Omega$, every point is written uniquely as $x=y+s\,\nu(y)$ with $y\in\partial\Omega$; the parallel hypersurfaces $\Gamma_s=\{s=\mathrm{const}\}$ have principal curvatures $\kappa_i/(1+s\kappa_i)$, so that their mean curvature $H_s$ satisfies $\partial_sH_s|_{s=0}=-|I\!I|^2$.

\begin{lemma}[Integral identities]\label{lem:pohozaev}
\[
\int_\Omega u=0,\qquad \int_\Omega|\nabla u|^2=\lambda\int_\Omega u^2,\qquad \int_\Omega u^2=\frac N2\,c^2|\Omega|.
\]
\end{lemma}

\begin{proof}
The first identity follows by integrating the equation, since $\partial_\nu u=0$; the second by multiplying it by $u$. For the third, multiply the equation by $x\cdot\nabla u$ (Rellich--Pohozaev). By \eqref{eq:gradzero} every boundary term containing $\nabla u$ vanishes, and one obtains
$\bigl(\tfrac N2-1\bigr)\int_\Omega|\nabla u|^2=\lambda\bigl(\tfrac N2\int_\Omega u^2-\tfrac{c^2}2\int_{\partial\Omega}x\cdot\nu\bigr)$. Since $\int_{\partial\Omega}x\cdot\nu=N|\Omega|$, the second identity gives the claim.
\end{proof}

These identities hold for every Schiffer domain, counterexamples included; they do not characterize the ball.

\begin{lemma}[Normal derivatives]\label{lem:normal}
On $\partial\Omega$,
\[
\partial_\nu^2u=-\lambda c,\qquad \partial_\nu^3u=\lambda cH,\qquad \partial_\nu^4u=\lambda c\bigl(\lambda-2|I\!I|^2-H^2\bigr).
\]
In dimension $2$, $\partial_\nu^4u=\lambda c(\lambda-3\kappa^2)$.
\end{lemma}

\begin{proof}
Near $\partial\Omega$, $\Delta u=\partial_s^2u+H_s\partial_su+\Delta_{\Gamma_s}u$. On $\partial\Omega$, $\partial_su=0$ and $\Delta_\Gamma u=0$, hence $\partial_s^2u=-\lambda c$. The coefficients of $s^0,s^1,s^2$ in the expansion of $u$ do not depend on the boundary point, so $\Delta_{\Gamma_s}u=O(s^3)$. Differentiating the equation once and twice in $s$ and setting $s=0$ gives $\partial_s^3u+H\partial_s^2u=0$ and $\partial_s^4u+2(\partial_sH_s)\partial_s^2u+H\partial_s^3u+\lambda\partial_s^2u=0$, whence the formulas.
\end{proof}

For a ball of radius $R$ one finds $\partial_\nu^4u=\lambda c(\lambda-(N^2-1)/R^2)$, which can be checked directly on the radial solution.

\begin{corollary}[Liu \cite{Liu}]\label{cor:liu}
$\partial_\nu^3u$ is constant on $\partial\Omega$ if and only if $\Omega$ is a ball.
\end{corollary}

\begin{proof}
By Lemma~\ref{lem:normal} this is equivalent to $H$ being constant; apply Alexandrov's theorem.
\end{proof}

\begin{remark}
These local relations do not constrain $\Omega$: since $\partial\Omega$ is analytic, the Cauchy--Kovalevskaya theorem provides, near any analytic hypersurface, a local solution with Cauchy data $(c,0)$. The obstruction is global.
\end{remark}

\begin{proposition}\label{prop:killing}
For every Killing field $X$ of $\R^N$, the function $Xu$ is either zero or a Dirichlet eigenfunction of $\Omega$ for the eigenvalue $\lambda$.
\end{proposition}

\begin{proof}
$X$ commutes with $\Delta$, and $Xu=0$ on $\partial\Omega$ by \eqref{eq:gradzero}.
\end{proof}

For translations this is the classical fact that $\lambda$ is also a Dirichlet eigenvalue of $\Omega$. It yields an unconditional lower bound.

\begin{proposition}\label{prop:krahn}
Every Schiffer domain in $\R^N$ satisfies $kR^*\ge 2^{1/N}j_{N/2-1,1}$.
\end{proposition}

\begin{proof}
The functions $\partial_iu$ do not all vanish identically, and $\int_\Omega\partial_iu=c\int_{\partial\Omega}\nu_i=0$. A nonzero $\partial_iu$ is therefore a sign-changing Dirichlet eigenfunction for $\lambda$ (Proposition~\ref{prop:killing}), so $\lambda$ is not the first Dirichlet eigenvalue of the domain $\Omega$ and $\lambda\ge\lambda_2^D(\Omega)$. By the Krahn--Szeg\H{o} inequality (see \cite[\S4.1]{Henrot}), $\lambda_2^D(\Omega)$ is at least the first Dirichlet eigenvalue of a ball of volume $|\Omega|/2$, whose radius is $2^{-1/N}R^*$; this gives the claim.
\end{proof}

\section{Shape derivatives and the boundary functional}\label{sec:shape}

Let $V$ be a deformation field, $\Omega_t=(I+tV)(\Omega)$, $\alpha=V\cdot\nu$ on $\partial\Omega$, and $u_t$ the normalized Neumann eigenfunction of $\Omega_t$, the eigenvalue being simple. We denote by $\dot u$ the material derivative and by $u'$ the shape derivative, so that $\dot u=u'+\nabla u\cdot V$. For the Dirichlet problem the condition $u_t=0$ on $\partial\Omega_t$ is transported by the deformation, so that $\dot u=0$ on the boundary; for the Neumann problem it is $\partial_{\nu_t}u_t=0$ that is transported, not the value of $u_t$. In \cite{Seck24} the functional
\[
J_2(\Omega)=\int_{\partial\Omega}u^2\dd\sigma+\gamma|\Omega|
\]
was studied under the assumption $\dot u=0$ on $\partial\Omega$. This assumption fails.

\begin{proposition}\label{prop:material}
Let $\Omega$ be a Schiffer domain, the eigenvalue being simple and $\|u\|_{L^2}=1$. It is not true that $\dot u=0$ on $\partial\Omega$ for every deformation preserving the volume to first order.
\end{proposition}

\begin{proof}
By \eqref{eq:gradzero}, $\dot u=u'$ on $\partial\Omega$. The shape derivative satisfies
\[
\partial_\nu u'=\dv_\Gamma(\alpha\nabla_\Gamma u)+\lambda u\,\alpha=\lambda c\,\alpha,\qquad
\lambda'=\int_{\partial\Omega}\bigl(|\nabla_\Gamma u|^2-\lambda u^2\bigr)\alpha=-\lambda c^2\int_{\partial\Omega}\alpha .
\]
Let $w$ solve $\Delta w+\lambda w=0$ in $\Omega$. Since $-\Delta u'-\lambda u'=\lambda'u$, Green's formula gives $-\lambda'\int_\Omega uw=\int_{\partial\Omega}(\partial_\nu u'\,w-u'\,\partial_\nu w)$. If $u'=0$ on $\partial\Omega$ whenever $\int_{\partial\Omega}\alpha=0$, then $\lambda'=0$ and $\int_{\partial\Omega}\alpha w=0$ for every $\alpha$ with zero mean, so $w$ is constant on $\partial\Omega$. Taking $w=e^{ik\omega\cdot x}$, the function $\omega\cdot x$ would take values in a discrete set on the connected boundary, hence would be constant on $\partial\Omega$ for every direction $\omega$, which is absurd.
\end{proof}

At a Schiffer domain, $dJ_2(\Omega)[V]=\int_{\partial\Omega}\bigl(2c\,u'+\alpha(Hc^2+\gamma)\bigr)$, and the term $2cu'$ cannot be dropped. If it were, criticality would give $H$ constant, hence a ball: the assumption $\dot u=0$ already contains the conclusion. The correct condition involves an adjoint state.

\begin{proposition}[Optimality condition]\label{prop:adjoint}
Let $u$ be the normalized eigenfunction ($\int_\Omega u^2=1$), $\lambda$ simple, and $\beta=\int_{\partial\Omega}u^2$; this choice of $\beta$ ensures that the Fredholm compatibility condition holds for the adjoint state $p$ defined by
\[
-\Delta p-\lambda p=-\beta u\ \text{in }\Omega,\qquad \partial_\nu p=u\ \text{on }\partial\Omega,\qquad \int_\Omega pu=0.
\]
Then $dJ_2(\Omega)[\alpha]=\int_{\partial\Omega}\alpha\bigl(Hu^2+\gamma-\beta u^2+2\lambda up-2\nabla_\Gamma u\cdot\nabla_\Gamma p\bigr)$. If $\Omega$ is a Schiffer domain, then $p=-\frac{\beta}{2\lambda}\,x\cdot\nabla u+c\,q$, where $q$, defined up to the addition of a multiple of $u$, solves
\[
\Delta q+\lambda q=0\ \text{in }\Omega,\qquad \partial_\nu q=1-\frac{|\partial\Omega|}{N|\Omega|}\,x\cdot\nu\ \text{on }\partial\Omega,
\]
and $\Omega$ is a critical point of $J_2$ (for some $\gamma$) if and only if $T:=H+2\lambda q$ is constant on $\partial\Omega$.
\end{proposition}

\begin{proof}
The general formula follows from Hadamard's formula and Green's formula applied to the pair $(u',p)$, using $2\int_\Omega uu'=-\int_{\partial\Omega}u^2\alpha$. At a Schiffer domain, $(-\Delta-\lambda)(x\cdot\nabla u)=2\lambda u$ and $D^2u=-\lambda c\,\nu\otimes\nu$ on $\partial\Omega$, which gives $\partial_\nu(x\cdot\nabla u)=-\lambda c\,x\cdot\nu$. The identities $\beta=c^2|\partial\Omega|$ and $c^2|\Omega|=2/N$ (Lemma~\ref{lem:pohozaev}) give the Neumann condition for $q$, whose compatibility follows from $\int_{\partial\Omega}x\cdot\nu=N|\Omega|$. Finally $u=c$, $\nabla_\Gamma u=0$ and $p=cq$ on $\partial\Omega$; adding a multiple of $u$ to $q$ changes $T$ by a constant.
\end{proof}

\section{The P-function}\label{sec:P}

Set $P=|\nabla u|^2+\lambda u^2$. The maximum principle for such P-functions goes back to Payne; see Sperb \cite{Sperb}. We record the precise form we need.

\begin{proposition}\label{prop:Pfunction}
On the open set $\{\nabla u\ne0\}$, let $B=(2\lambda u\nabla u-\nabla P)/|\nabla u|^2$ and $\kappa_N=\frac{N-2}{2(N-1)}$. Then
\[
\Delta P+B\cdot\nabla P+\kappa_N\frac{(\nabla u\cdot\nabla P)^2}{|\nabla u|^4}\ge0,
\]
with equality at a point if and only if the restriction of $D^2u$ to $(\nabla u)^\perp$ is a multiple of the identity there. In particular, for $N=2$,
\begin{equation}\label{eq:Pidentity}
\Delta P+\frac{2\lambda u}{|\nabla u|^2}\nabla u\cdot\nabla P-\frac{|\nabla P|^2}{|\nabla u|^2}=0 .
\end{equation}
\end{proposition}

\begin{proof}
By Bochner's formula and the equation, $\Delta P=2|D^2u|^2-2\lambda^2u^2$. Let $\rho=|\nabla u|$, $e_1=\nabla u/\rho$, and write in an adapted orthonormal frame $D^2u=\begin{pmatrix}a&b^T\\ b&C\end{pmatrix}$, $b\in\R^{N-1}$, $C\in\mathrm{Sym}(N-1)$, with $a+\operatorname{tr}C=-\lambda u$. Setting $A=a+\lambda u$ we get $\operatorname{tr}C=-A$, $\nabla P=2\rho(Ae_1+b)$, hence $\nabla u\cdot\nabla P=2\rho^2A$ and $|\nabla P|^2=4\rho^2(A^2+|b|^2)$. By Cauchy--Schwarz, $|C|^2\ge A^2/(N-1)$, with equality if and only if $C$ is scalar. Therefore
\[
\Delta P=2\bigl(a^2+2|b|^2+|C|^2\bigr)-2\lambda^2u^2\ge 4A^2+4|b|^2-4\lambda uA-2A^2\frac{N-2}{N-1}.
\]
Finally $4A^2+4|b|^2-4\lambda uA=-B\cdot\nabla P$ and $2A^2\frac{N-2}{N-1}=\kappa_N(\nabla u\cdot\nabla P)^2/\rho^4$. For $N=2$, $\kappa_2=0$ and $C$ is a scalar, which gives \eqref{eq:Pidentity}.
\end{proof}

The quadratic term equals $\kappa_N\rho^{-4}(\nabla u\cdot\nabla P)\,\nabla u\cdot\nabla P$, so $P$ is a subsolution of the linear operator $\Delta+\tilde B\cdot\nabla$ with $\tilde B=B+\kappa_N\rho^{-4}(\nabla u\cdot\nabla P)\nabla u$, locally bounded on each component of $\{\nabla u\ne0\}$. Recall that $\partial\Omega$ belongs to the critical set of $u$.

\begin{corollary}\label{cor:maxP}
The maximum of $P$ on $\overline\Omega$ is attained on $\partial\Omega$, where $P=\lambda c^2$, or at a critical point of $u$, where $P=\lambda u^2$. In dimension $2$, the same holds for the minimum of $P$.
\end{corollary}

\begin{proof}
If the maximum is attained at a point $x$ where $\nabla u(x)\ne0$, the strong maximum principle gives $P\equiv P(x)$ on the component of $\{\nabla u\ne0\}$ containing $x$, whose boundary consists of critical points of $u$ (possibly on $\partial\Omega$). In dimension $2$, \eqref{eq:Pidentity} is an equation, and the argument applies to the minimum.
\end{proof}

\begin{lemma}\label{lem:Pboundary}
Near $\partial\Omega$, with $s<0$ inside, $P=\lambda c^2-\tfrac23\lambda^2c^2H\,s^3+O(s^4)$. Just inside $\partial\Omega$, $P-\lambda c^2$ has the sign of $H$.
\end{lemma}

\begin{proof}
By Lemma~\ref{lem:normal}, $u=c-\frac{\lambda c}{2}s^2+\frac{\lambda cH}{6}s^3+O(s^4)$, and the tangential gradient of $u$ is $O(s^3)$.
\end{proof}

\begin{proposition}\label{prop:critical}
Let $\Omega\subset\R^N$ be a Schiffer domain.
\begin{enumerate}[label=(\roman*)]
\item There is a critical point $x_0\in\Omega$ of $u$ with $|u(x_0)|>|c|$.
\item If $N=2$ and $\Omega$ is not convex, there is a critical point $x_1\in\Omega$ of $u$ with $|u(x_1)|<|c|$.
\end{enumerate}
\end{proposition}

\begin{proof}
(i) At a point of $\partial\Omega$ farthest from the origin all principal curvatures are positive, so $H>0$ there. By Lemma~\ref{lem:Pboundary}, $\max P>\lambda c^2$, and this maximum is attained inside. If it is attained at a non-critical point, $P$ is constant on the corresponding component of $\{\nabla u\ne0\}$; this component does not touch $\partial\Omega$, where $P=\lambda c^2$, so its boundary contains a critical point $x_0\in\Omega$ with $\lambda u(x_0)^2=\max P>\lambda c^2$.
(ii) A non-convex planar domain has a boundary point where $\kappa<0$; the same argument, applied to the minimum of $P$ (Corollary~\ref{cor:maxP}), gives the claim.
\end{proof}

\begin{remark}
In dimension $2$, Proposition~\ref{prop:critical} shows that a single critical point forces convexity. In fact, a single critical point forces the disc, by the theorem of Willms and Gladwell \cite{WillmsGladwell} on saddle points; see Corollary~\ref{cor:single}. For $N\ge3$ the minimum principle is lost, and we do not know whether a non-convex Schiffer domain forces a critical point where $|u|<|c|$.
\end{remark}

\begin{theorem}\label{thm:bounds}
Let $\Omega\subset\R^N$ be a Schiffer domain and $M=\max_{\overline\Omega}|u|$. Then
\begin{enumerate}[label=(\alph*)]
\item $M>\sqrt N\,|c|$;
\item $|\nabla u|^2+\lambda u^2\le\lambda M^2$ in $\Omega$;
\item if $|u(x^*)|=M$, then $\operatorname{dist}(x^*,\partial\Omega)\ge\bigl(\tfrac\pi2\mp\arcsin\tfrac{|c|}{M}\bigr)/k$, with the sign $-$ if $u(x^*)$ and $c$ have the same sign and $+$ otherwise.
\end{enumerate}
\end{theorem}

\begin{proof}
(b) By Corollary~\ref{cor:maxP}, $\max P\le\max(\lambda c^2,\lambda M^2)=\lambda M^2$.
(a) By Lemma~\ref{lem:pohozaev}, $\frac1{|\Omega|}\int_\Omega P=\frac{2\lambda}{|\Omega|}\int_\Omega u^2=N\lambda c^2$. This mean is at most $\lambda M^2$ by (b), and the inequality is strict because $P$ is not constant: it equals $\lambda c^2<N\lambda c^2$ on $\partial\Omega$.
(c) By (b), $w=\arcsin(u/M)$ satisfies $|\nabla w|\le k$ where $|u|<M$; integrate along the segment joining $x^*$ to a nearest boundary point, on which $w=\arcsin(c/M)$, while $w(x^*)=\pm\pi/2$.
\end{proof}

These statements hold for all Schiffer domains, counterexamples included: they are necessary conditions and do not characterize the ball.

\section{A Talenti-type rigidity theorem}\label{sec:talenti}

The Schiffer balls are exactly the balls $B_R$ with $J_{N/2}(kR)=0$, with the radial solution $u=\zeta(r)/\zeta(R)$, where
\[
\zeta(r)=r^{1-N/2}J_{N/2-1}(kr),\qquad \zeta'(r)=-k\,r^{1-N/2}J_{N/2}(kr).
\]

\begin{theorem}\label{thm:talenti}
Let $\Omega\subset\R^N$, $N\ge2$, be a bounded domain with $C^2$ boundary and $u\in C^2(\overline\Omega)$ a solution of \eqref{eq:schiffer} with $c>0$. Assume the sign condition
\[
\textup{(H1)}\qquad u<c\quad\text{in }\Omega .
\]
Then $kR^*\ge j_{N/2,1}$, with equality if and only if $\Omega$ is a ball.
\end{theorem}

If $c<0$, the theorem applies to $-u$ under the condition $u>c$. Note that $\Omega$ is not assumed to be simply connected.

\begin{proof}
We normalize $c=1$.

\emph{Step 0.} By Lemma~\ref{lem:pohozaev}, $\int_\Omega u=0$, so $m_0:=\min u<0$. Since $u$ is analytic and nonconstant, its critical set and each of its level sets have measure zero.

\emph{Step 1 (role of (H1)).} For $t\in(m_0,1)$ let $E_t=\{u<t\}$. By (H1) and since $u=1$ on $\partial\Omega$, the set $\{u\le t\}$ is a compact subset of $\Omega$, so $\partial E_t=\{u=t\}$; moreover $\bigcup_{t<1}E_t=\Omega$. For almost every $t$ (Sard), $\{u=t\}$ is a smooth compact hypersurface with outward normal $\nabla u/|\nabla u|$. Without (H1), $\partial E_t$ would contain $\partial\Omega$ for $t>1$ and the argument below would fail.

\emph{Step 2 (flux).} For regular $t$, $F(t):=\int_{\{u=t\}}|\nabla u|\dd\Hm=\int_{E_t}\Delta u=-\lambda\int_{E_t}u>0$. Indeed, if $t\le0$ then $\int_{E_t}u<0$; if $t>0$ then $\int_{E_t}u=-\int_{\{u\ge t\}}u<0$.

\emph{Step 3 (distribution function).} By the coarea formula, $\mu(t)=|E_t|$ is absolutely continuous and strictly increasing, with $\mu'(t)=\int_{\{u=t\}}|\nabla u|^{-1}\dd\Hm>0$ for almost every $t$, $\mu(m_0)=0$ and $\mu(1^-)=|\Omega|$.

\emph{Step 4 (Cauchy--Schwarz and isoperimetry).} Let $L(t)=\Hm(\{u=t\})$. Then $L(t)^2\le\mu'(t)F(t)$, and the isoperimetric inequality $L\ge N\omega_N^{1/N}\mu^{1-1/N}$ gives, for almost every $t$,
\begin{equation}\label{eq:isop}
C_N\,\mu(t)^{2-2/N}\le L(t)^2\le\mu'(t)F(t)=-\lambda\,\mu'(t)\int_{E_t}u,\qquad C_N=N^2\omega_N^{2/N}.
\end{equation}

\emph{Step 5 (rearrangement).} Let $\varphi:(0,|\Omega|)\to(m_0,1)$ be the inverse of $\mu$. Since $\mu'>0$ almost everywhere, $\varphi$ is absolutely continuous, with $\varphi'(a)=1/\mu'(\varphi(a))$ almost everywhere; it extends continuously to $[0,|\Omega|]$ with $\varphi(0)=m_0$ and $\varphi(|\Omega|)=1$. The function $u^\sharp(x)=\varphi(\omega_N|x|^N)$ is the increasing radial rearrangement of $u$. With $\Phi(a)=\int_0^a\varphi$ we have $\int_{E_t}u=\Phi(\mu(t))$, and \eqref{eq:isop} becomes
\begin{equation}\label{eq:psi}
\psi(a):=C_N\,a^{2-2/N}\varphi'(a)+\lambda\Phi(a)\le0\quad\text{for almost every }a\in(0,|\Omega|),
\end{equation}
with $\Phi(0)=0$ and $\Phi(|\Omega|)=\int_\Omega u=0$. In the variable $a=\omega_Nr^N$ the radial Laplacian reads $\frac{d}{da}\bigl(C_Na^{2-2/N}\frac{d}{da}\bigr)$, so that $\psi'=\Delta u^\sharp+\lambda u^\sharp$, and $\psi\equiv0$ expresses that $u^\sharp$ solves the radial equation.

\emph{Step 6 (test function).} Let $Z(a)=\zeta(r)$ with $a=\omega_Nr^N$. Then $\bigl(C_Na^{2-2/N}Z'\bigr)'=-\lambda Z$ and
\[
C_Na^{2-2/N}Z'(a)=N\omega_Nr^{N-1}\zeta'(r)=-N\omega_N\,k\,r^{N/2}J_{N/2}(kr),
\]
which tends to $0$ as $a\to0$.

\emph{Step 7 (exact identity).} Let $A=|\Omega|$. Two integrations by parts give
\[
\int_0^AC_Na^{2-2/N}\varphi'Z'=C_NA^{2-2/N}Z'(A)+\lambda\int_0^A\varphi Z,\qquad
\lambda\int_0^A\Phi Z'=\lambda\Phi(A)Z(A)-\lambda\int_0^A\varphi Z=-\lambda\int_0^A\varphi Z,
\]
hence
\begin{equation}\label{eq:identityI}
I:=\int_0^A\psi(a)Z'(a)\dd a=C_NA^{2-2/N}Z'(A)=-N\omega_N\,k\,(R^*)^{N/2}J_{N/2}(kR^*).
\end{equation}

\emph{Step 8 (sign).} Assume $kR^*\le j_{N/2,1}$. For $0<a<A$ we have $kr<j_{N/2,1}$, so $Z'(a)<0$ and, by \eqref{eq:psi}, $I\ge0$. But $J_{N/2}(kR^*)\ge0$, so $I\le0$ by \eqref{eq:identityI}. Thus $I=0$, which gives $J_{N/2}(kR^*)=0$, that is $kR^*=j_{N/2,1}$, and $\psi=0$ almost everywhere.

\emph{Step 9 (equality).} If $\psi=0$ almost everywhere, the inequalities in \eqref{eq:isop} are equalities for almost every $t$, so $E_t$ is a ball $B(c_t,r_t)$. These balls are nested, so $|c_t-c_s|\le r_t-r_s$ for $s<t$, and $r_t\to R^*$ as $t\to1$; the centres converge and $\Omega=\bigcup_{t<1}E_t$ is a ball.

\emph{Step 10 (the balls).} Conversely, let $\Omega=B_R$ with $kR=j_{N/2,n}$ and $u=\zeta/\zeta(R)$. The function $\zeta$ solves $(r^{N-1}\zeta')'+k^2r^{N-1}\zeta=0$, and $r^{N-1}\cdot k^2r^{N-1}$ is increasing; by the Sonine--P\'olya theorem \cite[Thm.~7.31.1]{Szego}, the absolute values of the successive extrema of $\zeta$ form a strictly decreasing sequence. These extrema are $r=0$ and $r=j_{N/2,i}/k$, and their signs alternate. For $n\ge2$, the extremum at $j_{N/2,n-2}/k$ (at $r=0$ if $n=2$) has the sign of $\zeta(R)$ and a larger absolute value, so (H1) fails. For $n=1$, $u$ is increasing in $r$ and (H1) holds. Hence equality characterizes the ball.
\end{proof}

\begin{remark}
The interval $(0,j_{N/2,1})$ is exactly where $Z$ is monotone; beyond it $Z'$ changes sign and the argument no longer concludes. The positivity of $J_{N/2}$ plays the role of the maximum principle in \cite{Seck24}. The approach is that of Talenti's comparison principle \cite{Talenti}; what is specific to the Schiffer problem is the vanishing of the boundary terms in \eqref{eq:identityI}, due to $\partial_\nu u=0$ and $\int_\Omega u=0$. For instance, on the planar Schiffer discs of order $n=2,3,4$ one has $\max u=3.33,\,1.61,\,4.58$.
\end{remark}

Table~\ref{tab:thresholds} compares the threshold of Theorem~\ref{thm:talenti} with the unconditional bound of Proposition~\ref{prop:krahn}.

\begin{table}[ht]
\centering
\begin{tabular}{@{}ccc@{}}
\toprule
$N$ & every Schiffer domain: $2^{1/N}j_{N/2-1,1}$ & under (H1): $j_{N/2,1}$\\
\midrule
$2$ & $3.401$ & $3.832$\\
$3$ & $3.958$ & $4.493$\\
$4$ & $4.557$ & $5.136$\\
$5$ & $5.162$ & $5.763$\\
$6$ & $5.765$ & $6.380$\\
\bottomrule
\end{tabular}
\caption{Lower bounds for $kR^*$: Proposition~\ref{prop:krahn} versus Theorem~\ref{thm:talenti}.}
\label{tab:thresholds}
\end{table}

\begin{corollary}\label{cor:single}
Assume that $u$ has a single critical point $x_0$ in $\Omega$.
\begin{enumerate}[label=(\roman*)]
\item If $N=2$, then $\Omega$ is a disc.
\item If $N\ge3$, then either $\Omega$ is a ball or $kR^*>j_{N/2,1}$.
\end{enumerate}
\end{corollary}

\begin{proof}
Assume $c>0$. Since $\int_\Omega u=0$ and $u=c$ on $\partial\Omega$, the function $u$ attains a negative interior minimum, which must be $x_0$. Near the boundary $u=c-\frac{\lambda c}2s^2+O(s^3)<c$; if $u\ge c$ somewhere, $u$ would also have an interior maximum, a second critical point. Hence (H1) holds, and Theorem~\ref{thm:talenti} gives (ii). In dimension $2$, $x_0$ is an extremum, so $u$ has no saddle point and $\Omega$ is a disc by the theorem of Willms and Gladwell \cite{WillmsGladwell}.
\end{proof}

\begin{problem}\label{prob:H1}
Does \textup{(H1)} imply that $\Omega$ is a ball? For $N\ge3$, does a single critical point force $kR^*\le j_{N/2,1}$, hence the ball?
\end{problem}

\section{Infinitesimal rigidity of the ball}\label{sec:infinitesimal}

\begin{lemma}\label{lem:transcendence}
Let $\nu\ge0$ be an integer or a half-integer. Every positive zero of $J_\nu$ is transcendental.
\end{lemma}

\begin{proof}
For integer $\nu$ this follows from Siegel's theorem \cite{Siegel}; see \cite[\S15.28]{Watson}. Let $\nu=n+\frac12$. Then $J_{n+1/2}(x)=\sqrt{2/(\pi x)}\,\bigl(p_n(1/x)\sin x-q_n(1/x)\cos x\bigr)$, where $p_n,q_n$ are polynomials with rational coefficients without common real zero (the spherical Bessel functions $j_n$ and $y_n$ have no common zero). Let $x>0$ be an algebraic zero. If $p_n(1/x)\ne0$, then $\tan x$ is algebraic, hence so is $e^{2ix}$. If $p_n(1/x)=0$, then $\cos x=0$ and $x\in\frac\pi2+\pi\mathbb Z$ is transcendental. Both cases contradict the Lindemann--Weierstrass theorem, by which $e^{2ix}$ is transcendental for algebraic $x\ne0$.
\end{proof}

We use Lommel's relation \cite[\S9.6]{Watson}: for $m\ge1$,
\begin{equation}\label{eq:lommel}
J_{\nu+m}(x)=R_{m-1,\nu+1}(x)\,J_{\nu+1}(x)-R_{m-2,\nu+2}(x)\,J_\nu(x),
\end{equation}
where $R_{-1,\cdot}=0$ and $R_{m,\mu}$ is a polynomial in $1/x$ of degree $m$, with rational coefficients when $\mu$ is rational and leading coefficient $2^m\Gamma(\mu+m)/\Gamma(\mu)\ne0$.

\begin{lemma}\label{lem:bourget}
Let $\nu\ge0$ be an integer or a half-integer and $m\ge1$ an integer. Then $J_\nu$ and $J_{\nu+m}$ have no common positive zero.
\end{lemma}

\begin{proof}
For integer $\nu$ this is Bourget's hypothesis, proved by Siegel \cite{Siegel} (see \cite[\S15.28]{Watson}). In general, let $x>0$ be a common zero. Positive zeros of $J_\nu$ are simple and $J_{\nu+1}=\frac\nu xJ_\nu-J_\nu'$, so $J_{\nu+1}(x)\ne0$, and \eqref{eq:lommel} gives $R_{m-1,\nu+1}(x)=0$. Thus $x$ is algebraic, contradicting Lemma~\ref{lem:transcendence}.
\end{proof}

\begin{theorem}\label{thm:infschiffer}
Let $(\Omega_\varepsilon)$ be a family of Schiffer domains, $C^1$ in $\varepsilon$ with values in $C^{2,\alpha}$, with $\Omega_0=B_R$ and $J_{N/2}(kR)=0$. Then the first-order normal displacement is a spherical harmonic of degree $0$ or $1$, that is, a dilation or a translation.
\end{theorem}

\begin{proof}
Write $\partial\Omega_\varepsilon=\{r=R(1+\varepsilon h)\}$ and $u_\varepsilon=u_0+\varepsilon v+o(\varepsilon)$, with $u_0=\zeta/\zeta(R)$, so that $u_0'(R)=0$ and $u_0''(R)=-k^2$. To first order, the conditions $u_\varepsilon=c_\varepsilon$ and $\partial_\nu u_\varepsilon=0$ give $v(R,\cdot)=\mathrm{const}$ and $\partial_rv(R,\cdot)=Rk^2h$; the variations of $k$ and $c$ only affect the radial mode. For $h=Y_\ell$, a spherical harmonic of degree $\ell\ge2$, we get $v=A\,r^{1-N/2}J_{\ell+N/2-1}(kr)\,Y_\ell$, hence $AJ_{\ell+N/2-1}(kR)=0$. Lemma~\ref{lem:bourget}, applied with $\nu=N/2$ and $m=\ell-1$, gives $J_{\ell+N/2-1}(kR)\ne0$. Thus $A=0$, then $\partial_rv=0\ne Rk^2Y_\ell$, a contradiction.
\end{proof}

This infinitesimal statement is consistent with the local rigidity results of Agranovsky \cite{Agranovsky} in the plane and of Kobayashi \cite{Kobayashi} and Canuto \cite{Canuto} in every dimension; Lemma~\ref{lem:transcendence} makes the case of odd $N$ explicit. Accordingly, the counterexamples do not branch off Schiffer balls: in the plane the counterexample of \cite{ColbrookStepaniants} branches off Wheeler's constant-flux family \cite{Wheeler}, and Guo's domains appear near Bessel quasi-resonances \cite{Guo}.

\begin{proposition}[Berenstein problem]\label{prop:berenstein}
Let $\nu=N/2-1$ and consider the unit ball with $k=j_{\nu,n}$, $n\ge1$, which carries the radial solution of \eqref{eq:berenstein}. A deformation $r=1+\varepsilon Y_\ell$, with $Y_\ell$ a spherical harmonic of degree $\ell\ge2$, is admissible to first order for \eqref{eq:berenstein} if and only if
\[
kJ'_{\ell+\nu}(k)+\frac N2\,J_{\ell+\nu}(k)=0 .
\]
This condition never holds. Hence the ball is infinitesimally rigid for \eqref{eq:berenstein} in every dimension.
\end{proposition}

\begin{proof}
Let $u_0$ be the radial solution, so that $u_0(1)=0$, $u_0'(1)=1$ and, by the radial equation, $u_0''(1)=-(N-1)$. To first order, the conditions $u_\varepsilon=0$ and $\partial_\nu u_\varepsilon=1$ on $r=1+\varepsilon h$ give $v(1)=-h$ and $\partial_rv(1)=(N-1)h$; the variation of $k$ only affects the radial mode. With $h=Y_\ell$ and $v=A\,g(r)Y_\ell$, $g(r)=r^{-\nu}J_{\ell+\nu}(kr)$, we need $Ag(1)=-1$ and $Ag'(1)=N-1$, which is possible if and only if $g'(1)+(N-1)g(1)=0$, that is $kJ'_{\ell+\nu}(k)+\frac N2J_{\ell+\nu}(k)=0$.

Since $J'_\mu=J_{\mu-1}-\frac\mu xJ_\mu$, the left-hand side equals $xJ_{\nu+\ell-1}(x)-(\ell-1)J_{\nu+\ell}(x)$ at $x=k$. By \eqref{eq:lommel}, it is $S(x)J_{\nu+1}(x)+T(x)J_\nu(x)$ with
\[
S(x)=x\,R_{\ell-2,\nu+1}(x)-(\ell-1)\,R_{\ell-1,\nu+1}(x),
\]
a polynomial in $x$ and $1/x$ with rational coefficients. $S$ is not identically zero: the first term contains powers of $1/x$ of degree at most $\ell-3$, so the coefficient of $x^{-(\ell-1)}$ is $-(\ell-1)\,2^{\ell-1}\Gamma(\nu+\ell)/\Gamma(\nu+1)\ne0$. At $x=j_{\nu,n}$ we have $J_\nu=0$ and $J_{\nu+1}\ne0$, so the condition is equivalent to $S(j_{\nu,n})=0$, which would make $j_{\nu,n}$ algebraic. This contradicts Lemma~\ref{lem:transcendence}.
\end{proof}

For $N=2$ the condition reads $kJ'_\ell(k)+J_\ell(k)=0$. Numerically, the relative value of the left-hand side stays above $10^{-3}$ for $N\le5$, $n\le100$ and $2\le\ell\le200$, in agreement with Proposition~\ref{prop:berenstein}. The counterexample of \cite{ColbrookSadeghiStepaniants} therefore does not branch off a disc either.

\section{A numerical illustration: Guo's domain in dimension four}\label{sec:numerics}

In dimension $N\ge3$, symmetric domains have non-isolated critical sets: if $\Omega$ is invariant under $O(N-1)$, a critical point off the axis generates a critical sphere $S^{N-2}$, and the Poincar\'e--Hopf formula becomes a Morse--Bott formula.

\begin{lemma}\label{lem:morsebott}
Assume that the critical set of $u$ in $\Omega$ is a finite union of critical submanifolds $C_i$, non-degenerate in the sense of Bott, with normal indices $\iota_i$. Then $\sum_i(-1)^{\iota_i}\chi(C_i)=1$.
\end{lemma}

\begin{proof}
Assume $c>0$. Near $\partial\Omega$, $\partial_su=-\lambda cs+O(s^2)>0$ for $s<0$: on an inner parallel hypersurface $\Gamma_{-\varepsilon}$ the field $\nabla u$ points outward. The region bounded by $\Gamma_{-\varepsilon}$ has the Euler characteristic of a ball (Alexander duality), and the Morse--Bott form of the Poincar\'e--Hopf formula gives the result.
\end{proof}

\subsection*{Planar reduction and method} Following \cite{Guo}, we look for $\Omega=\{(x_1,x')\in\R\times\R^3:(x_1,|x'|)\in D\}$, where $D=\psi(\mathbb D)$, $\mathbb D$ is the unit disc and $\psi(w)=\sum_{j\ \mathrm{odd}}c_jw^j$ has real coefficients; we normalize $\lambda=1$. For $u=u(x_1,y)$ with $y=|x'|$, one has $\Delta_{\R^4}u=u_{x_1x_1}+u_{yy}+\frac2yu_y$, and $F=yu$ satisfies $\Delta_{\R^2}F=y\,\Delta_{\R^4}u$. Problem \eqref{eq:schiffer} becomes the planar problem
\[
\Delta F+F=0\ \text{in }D,\qquad F=y\ \text{and}\ \nabla F=\nabla y\ \text{on }\partial D.
\]
On the disc, $V=F\circ\psi$ satisfies $\Delta V+|\psi'|^2V=0$, with $V=\operatorname{Im}\psi$ and $\partial_rV=\operatorname{Im}(w\psi'(w))$ for $|w|=1$. The symmetry imposes $V=\sum_{n\ \mathrm{odd}}v_n(r)\sin n\theta$. This sector excludes the Dirichlet eigenfunctions of Proposition~\ref{prop:killing}: $\partial_{x_1}u$ is odd in $x_1$, and the rotations mixing $x_1$ and $x'$ break the $O(3)$ invariance. We discretize by Chebyshev collocation in $r$, with the odd extension $V(-r,\theta)=-V(r,\theta)$, and by odd sines in $\theta$. The unknowns are $V$ and $c_1,c_3,\dots,c_{2h-1}$; the square system is solved by Newton's method, starting only from the values published in \cite[\S1.5]{Guo}, namely $c_1\approx14.840043$, $c_3\approx0.0912$, $c_5\approx0.0942$, $c_7\approx0.0964$, the other coefficients being set to zero. Newton's method converges quadratically in three iterations. These computations are not certified.

\begin{figure}[!htbp]
\centering
\includegraphics[width=0.55\textwidth]{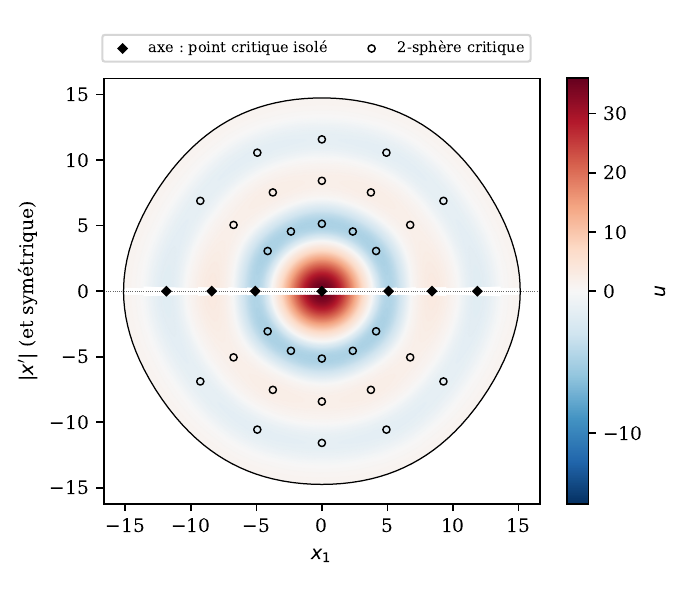}
\caption{Meridian section of Guo's domain in dimension $4$: the function $u$ in the plane $(x_1,|x'|)$, extended by symmetry ($c=1$, $\lambda=1$). Diamonds: isolated critical points on the axis; circles: critical $2$-spheres, each represented by its two traces in the plane.}
\label{fig:guo}
\end{figure}

\begin{table}[ht]
\centering
\small
\begin{tabular}{@{}lll@{}}
\toprule
Quantity & Value & Reference\\
\midrule
$c_1$ (stable from $(N,h)=(61,24)$ to $(121,48)$) & $14.8400430097$ & \cite{Guo}: $14.840043$\\
$c_3;\ c_5;\ c_7$ & $0.0911958214;\ 0.0942290934;\ 0.0963587384$ & \cite{Guo}: $0.0912;\ 0.0942;\ 0.0964$\\
Residual of the equations & $\le5\cdot10^{-8}$ & \\
Radius of $\partial D$ & $14.7213$ to $15.1343$ & \\
$\min_{|w|=1}\operatorname{Re}\bigl(1+w\psi''/\psi'\bigr)$ & $0.6958>0$ (convex) & \\
$\int_\Omega u$;\ $\int_\Omega u^2/|\Omega|$ & $2\cdot10^{-8}$ ($|\Omega|\approx2.37\cdot10^5$);\ $2.00000000$ & Lemma~\ref{lem:pohozaev}\\
$M=u(0)$ & $35.778$; (H1) fails & Thm.~\ref{thm:bounds}(a): $M>2$\\
$\max P/(\lambda M^2)$ away from the axis & $0.9736$ & Thm.~\ref{thm:bounds}(b)\\
Critical points on the axis & $7$: indices $4\,(\times3)$, $0\,(\times2)$, $3\,(\times2)$ & \\
Critical $2$-spheres & $15$: normal indices $0\,(\times5)$, $1\,(\times8)$, $2\,(\times2)$ & \\
Morse--Bott sum & $3+2-2+2(5-8+2)=1$ & Lemma~\ref{lem:morsebott}\\
\bottomrule
\end{tabular}
\caption{Guo's domain in dimension $4$, reconstructed independently ($c=1$, $\lambda=1$).}
\label{tab:guo}
\end{table}

\begin{table}[ht]
\centering
\small
\begin{tabular}{@{}ccccc@{}}
\toprule
$\theta$ (in $\mathbb D$) & $H$ & $\partial_\nu^3u$ & $1-2|I\!I|^2-H^2$ & $\partial_\nu^4u$\\
\midrule
$0.40$ & $0.20719$ & $0.20719$ & $0.927842$ & $0.927969$\\
$0.80$ & $0.19455$ & $0.19455$ & $0.936898$ & $0.937041$\\
$1.20$ & $0.21194$ & $0.21194$ & $0.925036$ & $0.925086$\\
$1.50$ & $0.19353$ & $0.19353$ & $0.937446$ & $0.937458$\\
\bottomrule
\end{tabular}
\caption{Normal derivatives on the boundary of the four-dimensional domain: Lemma~\ref{lem:normal} with $\lambda=c=1$. Here $H=\kappa+2\nu_y/y$ and $|I\!I|^2=\kappa^2+2(\nu_y/y)^2$, where $\kappa$ is the meridian curvature and $\nu_y/y$ the transverse curvature. The derivatives are obtained by polynomial fitting along the normal; the discrepancy of order $10^{-4}$ at order four comes from this fitting.}
\label{tab:normal}
\end{table}

\subsection*{Comments} The reconstruction is independent of the proof of \cite{Guo}: starting from four rounded coefficients, it recovers all published digits and refines them to ten digits that are stable under refinement of the discretization (the residual limits the guaranteed accuracy to about $10^{-7}$). The domain is convex and its frequency is high: $kR^*\approx14.8$, against $j_{2,1}\approx5.14$. It violates (H1), since $u(0)\approx35.8$, and therefore lies outside the scope of Theorem~\ref{thm:talenti}. Figure~\ref{fig:guo} shows the computed critical sets. The Morse--Bott sum equals $1$, as required by Lemma~\ref{lem:morsebott}; this is a consistency check of the computed critical sets, not a proof that none was missed, since pairs with cancelling contributions would go undetected.

\section{Open problems}\label{sec:open}

\begin{enumerate}[label=\arabic*.]
\item Does (H1) force the ball (Problem~\ref{prob:H1})? An upper bound for $kR^*$ under (H1) would suffice.
\item For $N\ge3$, does a non-convex Schiffer domain force a critical point where $|u|<|c|$? A substitute for the minimum principle would be needed.
\item Is there a threshold $m^*$ such that $\max|u|\le m^*|c|$ forces the ball? Theorem~\ref{thm:bounds} gives $\max|u|>\sqrt N|c|$ for every Schiffer domain.
\item Quantitative local rigidity of the ball among Schiffer domains with $kR^*\le K$, with an explicit neighbourhood.
\item Reconstruct Guo's domains in other dimensions and repeat the tests of Section~\ref{sec:numerics}.
\end{enumerate}

\subsection*{Code availability} The Python scripts used in Section~\ref{sec:numerics} (planar reduction, Newton's method, tests, Morse--Bott count) and for the Bessel computations of Section~\ref{sec:infinitesimal} are available from the author and will be deposited in a public repository.

\subsection*{Use of artificial intelligence} The numerical computations and part of the writing were carried out with the assistance of a language model (Claude, Anthropic). The author has reviewed the whole text and takes full responsibility for it.

\subsection*{Acknowledgements} The author thanks the colleagues of his laboratory for their encouragement.

\end{document}